\documentclass{crm-eca}

\usepackage{aut-free-groups}  
\usepackage{latexsym}
\usepackage{amssymb, amsfonts, amsmath}

\newtheorem{theorem}{Theorem}

\theoremstyle{definition}
\newtheorem{definition}[theorem]{Definition}

\theoremstyle{remark}

\begin{document}

\papertitle{Orbit decidability, applications and variations}


\paperauthor{Enric Ventura}
\paperaddress{Universitat Polit\`{e}cnica de Catalunya} 
\paperemail{enric.ventura@upc.edu}

\paperthanks{The author thanks the CRM hospitality during the research programme on Automorphisms of Free Groups. He
also acknowledges partial support from the Spanish Government through grant number MTM2011-25955.}

\makepapertitle

\section{Orbit decidability}

In many areas of mathematics and in innumerable topics and situations, the notion of \emph{transformation} plays an
important role. If $X$ is the set or collection of objects we are interested in, a transformation of $X$ is usually
understood to be just a map $\alpha \colon X\to X$. And whenever the context highlights a certain collection of
``interesting" maps $A\subseteq \operatorname{Map}(X,X)$ (namely, endomorphisms or automorphisms of $X$ if $X$ is an
algebraic structure, continuous maps or isometries of $X$ if $X$ is a topological or a geometric object, etc), one
naturally has the notion of orbit: the \emph{$A$-orbit} of a point $x\in X$ is the set of all its $A$-images $xA=\{
x\alpha \mid \alpha \in A\} \subseteq X$. In all these situations, there is a problem which is usually crucial when
studying algorithmic aspects of many of the interesting problems one can formulate about the objects in $X$ and how do
they relate to each other under the transformations in $A$: orbit decidability.

\begin{definition}
Let $X$ be a set, and let $A\subseteq \operatorname{Map}(X,X)$ be a set of transformations. We say that $A$ is
\emph{orbit decidable} (\emph{OD} for short) if there is an algorithm which, given $x,y\in X$, decides whether $x\alpha
=y$ for some $\alpha \in A$. (Sometimes the algorithm is required to provide such an $\alpha$, if it exists.)
\end{definition}

There are 
lots of examples of very classical algorithmic problems which are of this kind. For example,
the conjugacy problem of a group $G$ is just the orbit decidability for the set of inner automorphisms
$A=\operatorname{Inn}(G)$ (and recall that the word problem of $G$ is a special subproblem). The classical Whitehead
algorithm for the free group $F_n$ is just a solution to the orbit decidability of the full automorphism group
$A=\operatorname{Aut}(F_n)$, and all the variations of this problem (replacing elements to conjugacy classes or
subgroups, of tuples of them, etc; replacing automorphisms to certain kind of automorphisms or endomorphisms, etc;
moving to other families of groups $G$ or algebraic structures, etc) are nothing else than other instances of orbit
decidability.

A recent result by Bogopolski-Martino-Ventura~\cite{BMV} gave a renovated protagonism to the notion of orbit
decidability. We first remind a couple of other concepts. The \emph{twisted conjugacy problem} (\emph{TCP}) for a group
$G$ consists on deciding, given $\alpha \in \operatorname{Aut}(G)$ and two elements $u,v\in G$, whether there exists
$x\in G$ such that $(x\alpha)^{-1}ux=v$; note that if $\alpha$ is the identity this is precisely the standard conjugacy
problem (\emph{CP}) for $G$ but, in general, it is a strictly stronger algorithmic problem
(see~\cite[Corollary~4.9]{BMV} for an example of a group with solvable \emph{CP} but unsolvable \emph{TCP}). On the
other hand, for a short exact sequence of groups, $1\longrightarrow F\stackrel{\alpha}{\longrightarrow}
G\stackrel{\beta}{\longrightarrow} H\longrightarrow 1$, and since $F\alpha$ is a normal subgroup of $G$, for every
$g\in G$, the conjugation $\gamma_g$ of $G$ induces an automorphism of $F$, $\varphi_g \colon F\to F$, $x\mapsto
g^{-1}xg$ (which does not necessarily belong to $Inn(F)$). The set of all such automorphisms, $A_{G} =\{ \varphi_g \mid
g\in G\}$, is a subgroup of $Aut(F)$ called the \emph{action subgroup} of the given short exact sequence.

\begin{theorem}[Bogopolski-Martino-Ventura~\cite{BMV}]\label{ses}
Let $1\longrightarrow F\stackrel{\alpha}{\longrightarrow} G\stackrel{\beta}{\longrightarrow} H\longrightarrow 1$ be a
short exact sequence of groups (given by finite presentations and the images of generators) such that
\begin{itemize}
\item[(i)] $F$ has solvable TCP,
\item[(ii)] $H$ has solvable CP, and
\item[(iii)] for every $1\neq h\in H$, the subgroup $\langle h\rangle$ has finite index in its centralizer $C_H(h)$, and
there is an algorithm which computes a finite set of coset representatives, $z_{h,1},\ldots ,z_{h,t_h}\in H$ (i.e.,
$C_H(h)=\langle h \rangle z_{h,1}\sqcup \cdots \sqcup \langle h \rangle z_{h,t_h}$).
\end{itemize}
Then,
 $$
G \text{ has solvable CP } \Longleftrightarrow A_G =\{ \varphi_g \mid g\in G\} \leqslant Aut(F) \text{ is OD.}
 $$
\end{theorem}

Hypothesis (iii) is somehow restrictive, but at the same time satisfied by many groups: for example, free groups (where
the centralizer of an element $1\neq h$ is cyclic and generated by its maximal root) and it is not difficult to see
that torsion-free hyperbolic groups also satisfy it, see~\cite[Subsection~4.2]{BMV}.

The correct way to think about this theorem is the following: it reduces the \emph{CP} for a group $G$ to the
\emph{TCP} plus a certain \emph{OD} problem for a certain subgroup $H\leqslant G$. It is true that the \emph{TCP} is
harder than the standard \emph{CP}, and the resulting \emph{OD} problem is sometimes more technical than the original
problem; but both of them take place \emph{in the subgroup $H$} rather than in $G$. In all situations when $H$ is a
group significantly easier than $G$, Theorem~\ref{ses} reduces the \emph{CP} for $G$ to two independent problems, maybe
more technical but in an easier group $H$. Let us say in a different way: for any group $H$ where one knows how to
solve the \emph{TCP}, Theorem~\ref{ses} gives a great tool to investigate the solvability/unsolvability of the
\emph{CP} in a vast family of extensions of $H$, by means of finding orbit decidable/orbit undecidable subgroups of
$\operatorname{Aut}(H)$.

\section{Applications}

The idea behind Theorem~\ref{ses} has proven to be quite fruitful, being the starting point of a collection of papers
and preprints. The first one was~\cite{BMMV}, where Bogopolski-Martino-Maslakova-Ventura solved $TCP(F_n)$; combining
this with Brinkmann's result that cyclic subgroups of $\operatorname{Aut}(F_n)$ are \emph{OD} (see \cite{Br}), one
immediately gets a solution to the \emph{CP} for free-by-cyclic groups. (We remark that all these arguments made a
crucial use of a result of Maslakova~\cite{Mas} on computability of the fixed subgroup of an automorphism of a free
group, which is now under revision because of incorrectness of the original argument, see~\cite{BM}; for an alternative
solution to the \emph{CP} for free-by-cyclic groups given by Bridson-Groves, see~\cite{BrGr}.)

In Theorem~\ref{ses}, we can take both $F$ and $H$ to be free groups. But a well known construction due to C. Miller,
see~\cite{Mi}, provided examples of free-by-free groups with unsolvable \emph{CP}. Hence, Theorem~\ref{ses} tells us
that $\operatorname{Aut}(F_n)$ must contain orbit undecidable subgroups $A\leqslant \operatorname{Aut}(F_n)$. This is
not the case in rank 2 (every finitely generated subgroup of $\operatorname{Aut}(F_2)$ is \emph{OD},
see~\cite[Proposition~6.13]{BMV}), but they certainly do exist for higher rank, $n\geqslant 3$. A closer look to these
negative examples revealed a general way to construct orbit undecidable subgroups inside $\operatorname{Aut}(G)$, as
soon as $F_2 \times F_2$ embeds into it (see~\cite[Section 7]{BMV}). This allowed to construct lots of new extensions
of groups with unsolvable conjugacy problem. For example, since $F_2$ embeds in $GL_2(\mathbb{Z})$, $F_2\times F_2$
embeds in $GL_4(\mathbb{Z})$ and one can deduce that $GL_4(\mathbb{Z})=\operatorname{Aut}(\mathbb{Z}^4)$ contains orbit
undecidable subgroups which, via Theorem~\ref{ses}, implies the existence of $\mathbb{Z}^n$-by-free groups with
$n\geqslant 4$ and unsolvable \emph{CP} (see~\cite[Proposition~7.5]{BMV}). At this point it is worth mentioning that
non of these arguments apply to the case of dimension 3 so, at the time of writing, it is an open problem whether there
exists $\mathbb{Z}^3$-by-free groups with unsolvable \emph{CP} (i.e. whether or not $GL_3(\mathbb{Z})$ contains orbit
undecidable subgroups).

These last results were used by Sunic-Ventura in~\cite{SV} to see that there exist automaton groups (i.e. subgroups of
the automorphism group of a regular rooted tree, generated by finite self-similar sets) with unsolvable \emph{CP}. In
fact, in~\cite{SV} and using techniques of Brunner and Sidki, it was proved that $\mathbb{Z}^d \rtimes \Gamma$ is an
automaton group for every finitely generated $\Gamma \leqslant GL_d(\mathbb{Z})$. Then, by modifying the construction
in~\cite{BMV} at the cost of increasing the dimension in 2 units, a finitely generated, orbit undecidable, free
subgroup $\Gamma$ of $GL_d(\mathbb{Z})$ was constructed, for $d\geqslant 6$. Using both results together with
Theorem~\ref{ses}, one gets automaton groups with unsolvable \emph{CP} (and additionally being [free-abelian]-by-free).

In the preprint~\cite{GmV}, Gonz\'{a}lez-Meneses and Ventura consider the braid group $B_n$ and solve \emph{TCP}($B_n$).
With a first superficial look, it may seem an easy problem because it is well known that $\operatorname{Out}(B_n)\simeq
C_2$, with the non-trivial element represented by the automorphism $\alpha \colon B_n \to B_n$ which inverts all
generators, $\sigma_i \mapsto \sigma_i^{-1}$. However, the conjugacy problem twisted by this $\alpha$ (namely solving
the equation $(x\alpha)^{-1}ux=v$ for $x\in B_n$) becomes a quite delicate combinatorial problem about palindromic
braids (see~\cite{GmV} for details). Furthermore, it is easy to see that every finitely generated subgroup $A\leqslant
\operatorname{Aut}(B_n)$ is orbit decidable; hence, every extension of $B_n$ by a torsion-free hyperbolic group $H$ has
solvable \emph{CP}, see~\cite[section~5]{GmV}.

A kind of opposite situation happens in Thompson's group $F$. Here, the automorphism group is quite big; but it is
known that every automorphism of $F$ can be realized as the conjugation by some element in $\widetilde{EP}_2$ (a
certain discrete subgroup of $\operatorname{Homeo}([0,1])$ containing $F$). So, $F$ has lots of automorphisms, but they
all are structurally easy. This allowed Burillo-Matucci-Ventura to solve $TCP(F)$ in~\cite{BuMV}. Since it is also
proved there that $F_2\times F_2$ does embed in Thomson's group $F$, one deduces the existence of Thompson-by-free
groups with unsolvable \emph{CP}.

A similar project is currently being carried over by Fern\'{a}ndez-Alcober, Ventura and Zugadi for the family of
Grigorchuk-Gupta-Sidki groups,~\cite{FVZ}.

We encourage the (algorithmic oriented) reader to push the same idea further into his own area of expertise: choose
your favorite group $G$, and try to solve \emph{TCP}($G$). This will not be a very interesting result by itself (it is
just a technical variation of \emph{CP}($G$)), but it will pave the way (via Theorem~\ref{ses}) to study the \emph{CP}
in a vast collection of extensions of $G$: you will have chances to prove results of the type ``all
$G$-by-[torsion-free hyperbolic] groups have solvable \emph{CP}", or ``there exists a $G$-by-free group with unsolvable
\emph{CP}".

\section{Variations on orbit decidability}

The definition of orbit decidability admits variations, pointing to deeper algorithmic problems. We present here one of
these possible variations that we find interesting. It is not totally clear, by the moment, whether is it related to
some algebraic problem, like standard orbit decidability is related to the \emph{CP} via Theorem~\ref{ses}. Even if it
is not, the problems it provides are interesting enough by themselves.

\begin{definition}
Let $G$ be a group, and $A\leqslant \operatorname{Aut}(G)$. We say that $A$ is ($m$-)\emph{subgroup orbit decidable},
\emph{($m$-)SOD} for short, if there is an algorithm which, given $g, h_1,\ldots ,h_m \in G$, decides whether $g\alpha
\in H=\langle h_1,\ldots ,h_m\rangle \leqslant G$ for some $\alpha \in A$.
\end{definition}

Since in $F_n$, as well as in $\mathbb{Z}^n$, roots of elements are well-defined and must be preserved by automorphisms
(i.e. $x\alpha =y$ implies $\hat{x}\alpha =\hat{y}$), it is easy to see that, for every $A$, solvability of
\emph{OD}(A) implies solvability of $1$-\emph{SOD}($A$). However, $m$-\emph{SOD}($A$) for $m\geqslant 2$ looks like a
much more complicated problem, even over the free abelian group.

Over the free group $F_n$, two special instances of this problem are solved in the literature. Silva-Weil solved
in~\cite{SW} the problem \emph{SOD}($\operatorname{Aut}(F_2)$): given an element $x$ and a subgroup $H$ of the rank two
free group $F_2$, one can algorithmically decide whether $x\alpha \in H$ for some $\alpha \in \operatorname{Aut}(F_2)$.
And Clifford-Goldstein \cite{CG} gave an algorithm solving the particular case of \emph{SOD}($\operatorname{Aut}(F_n)$)
where the given input $x$ is a primitive element: there is an algorithm deciding whether a given subgroup $H\leqslant
F_n$ contains a primitive element of $F_n$. The rest of the problem \emph{SOD}($\operatorname{Aut}(F_n)$) remains open,
and nothing is know for other subgroups $A\leqslant \operatorname{Aut}(F_n)$.


Over the free abelian group $\mathbb{Z}^n$, \emph{SOD}($GL_n(\mathbb{Z})$) is an exercise (just a matter of
$\operatorname{gcd}$'s of the entries of the involved vectors). But, for a fixed given matrix $A\in GL_n(\mathbb{Z})$,
the problem \emph{SOD}($\langle A\rangle$) is much more interesting: after projectivizing $\mathbb{Z}^n$, the
automorphism $A\colon \mathbb{Z}^n \to \mathbb{Z}^n$ induces a map $\varphi \colon \mathbb{P}^{n-1}(\mathbb{Z}) \to
\mathbb{P}^{n-1}(\mathbb{Z})$, and \emph{SOD}($\langle A\rangle$) becomes the problem of deciding whether a given orbit
of $\varphi$ intersects a given (projective) linear variety in $\mathbb{P}^{n-1}(\mathbb{Z})$ (for $n=2$, this problem
becomes a nice exercise in linear algebra, involving the eigenvalues of $A$).

\end{document}